\documentclass[11pt]{article}

\usepackage{amssymb}
\usepackage{color,amsmath,amssymb, amsfonts,amstext,amsthm}

\usepackage{epsfig, graphicx, graphics}
\usepackage{longtable}

\newcommand{\om}{\omega}

\renewcommand{\phi}{\varphi}

\renewcommand{\a}{\alpha}

\newcommand{\R}{{\mathbb R}}

\newcommand{\eps}{\varepsilon}

\newcommand{\EX}{{\mathbb{E}}}

\title{Mean exit time and escape probability for dynamical systems
driven by L\'evy noise
\footnote{This work was partly supported by the NSF Grants 0511411,
0620539 and 0923111, the Simons Foundation grant 208236, and the NSFC grants 10971225 and 11028102. } }

\author{Ting Gao$^{1}$, Jinqiao Duan$^{1, 2}$, Xiaofan Li$^{1}$, Renming Song$^{3}$    \\
\\
1.  Department of Applied Mathematics\\ Illinois Institute of Technology \\
  Chicago, IL 60616, USA \\
 \emph{E-mail: tinggao0716@gmail.com,  duan@iit.edu,  lix@iit.edu}\\
2. Institute for Pure and Applied Mathematics, University of California\\
Los Angeles, CA 90095, USA\\
\emph{E-mail: jduan@ipam.ucla.edu }\\
3. Department of Mathematics\\ University of Illinois at
Urbana-Champaign\\
 Urbana, IL 61801, USA\\
 \emph{E-mail: 	rsong@math.uiuc.edu} }

\begin{document}
\date\today

\maketitle

\pagestyle{plain}

\begin{abstract}
The mean first exit time and escape probability are utilized to quantify dynamical behaviors of stochastic differential equations with non-Gaussian $\alpha-$stable type L\'evy motions.  Both deterministic quantities are characterized by differential-integral equations (i.e., differential equations with nonlocal terms) but with different exterior conditions.
The non-Gaussianity of noises manifests as nonlocality at the level of mean exit time and escape probability. An objective of this paper   is to make mean exit time and escape probability as efficient computational tools, to the applied probability community, for quantifying
stochastic dynamics.
An accurate numerical scheme is developed and validated for computing the mean   exit time and escape probability.
Asymptotic solution for the mean exit time is given when the pure jump measure
in the L\'evy motion is small.


From both the analytical and numerical results,
it is observed that the mean exit
time depends strongly on the domain size
and the value of $\a$ in the $\alpha-$stable L\'evy jump measure.
The mean exit time can measure which of the two competing
factors in $\a$-stable L\'evy motion, i.e. the jump frequency or the jump size,
is dominant in helping a process exit a bounded domain.
The escape probability is shown to vary with the underlying vector field (i.e., drift).
The mean exit time and escape probability could become discontinuous
at the boundary of the domain, when the process is subject to certain
deterministic potential and the value of $\a$ is in $(0,1)$.






\medskip




\bigskip

\end{abstract}

\section{Motivation}  \label{intro}

Random fluctuations in  complex systems in engineering and science are often  non-Gaussian \cite{Woy, Diego1, Diego2}.
For instance, it has been argued that diffusion by
geophysical turbulence \cite{Shlesinger} corresponds, loosely
speaking, to a series of  ``pauses", when the particle is trapped by
a coherent structure, and ``flights" or ``jumps" or other extreme
events, when the particle moves in the jet flow. Paleoclimatic data
\cite{Dit} also indicate such irregular processes.

L\'evy motions are thought to be appropriate models for
non-Gaussian processes with jumps \cite{Sato-99}.
Recall that a L\'evy motion $L(t)$, or $L_t$, is a stochastic process with stationary
and independent increments. That is, for any $s, t$ with  $0\le s< t$, the distribution
of $L_t-L_s$ only depends on $t-s$, and for any $0\le t_0<t_1<\cdots<t_n$,
$L_{t_i}-L_{t_{i-1}}$, $i=1, \cdots, n$, are independent. Without loss of
generality, we may assume that the sample paths of $L_t$ are almost surely
right continuous with left limits.

This generalizes the Brownian motion $B(t)$, which satisfies all
these three conditions.  But \emph{additionally},  (i) almost all
sample paths of the Brownian motion     are  continuous in time in the
usual sense and (ii) Brownian motion's increments are Gaussian
distributed.

  Stochastic differential equations (SDEs)
with   non-Gaussian L\'evy noises
have attracted much attention  recently \cite{Apple,  Schertzer}.
To be specific, let us   consider the following scalar SDE
with a non-Gaussian L\'evy motion
\begin{equation} \label{sde}
{\rm d}X_{t}  =  f (X_{t})  {\rm d}t +  {\rm d} L_{t}, \;\; X_0 = x,
\end{equation}
 where $f$ is a vector field (or drift), and $L_{t}$ is a scalar L\'evy motion defined in a probability space $(\Omega, \mathcal{F}, \mathbb{P})$.

\vskip 12pt

 We   study the first
exit problem for the solution process $X_{t}$ from bounded domains.
The exit phenomenon, i.e., escaping from a bounded domain in state space,
 is an impact of randomness on the evolution of such dynamical systems.
 Two concepts are applied to quantify the exit phenomenon: \emph{mean exit time} and \emph{escape probability}.

 We define the  first exit time from  the spatial domain $D$  as follows:
\[
\tau (\omega):= \inf \{t \geq 0, X_{t}(\om, x)  \notin D \},
\]
and the mean exit time is denoted as $u(x) := \EX \tau$.
The  likelihood of a particle (or a solution path $X_t$), starting at a point $x$,
first escapes a domain $D$ and   lands
in a subset $E$ of $D^c$ (the complement of $D$) is called escape probability and is denoted as $P_E(x)$.

The existing work on mean exit time gives asymptotic estimate for $u(x)$
when the noise intensity is sufficiently small, i.e.,
the noise term in \eqref{sde} is $\eps \; {\rm d}L_t$ with $0<\eps \ll 1$.
See, for example, Imkeller and Pavlyukevich \cite{ImkellerP-06,
ImkellerP-08}, and Yang and Duan \cite{YangDuan}.

In the present paper, however, we numerically investigate
mean exit time and escape probability for arbitrary noise intensity.
 The mean exit time $u(x)$ and escape probability $P_E(x)$ for a dynamical system, driven by a
non-Gaussian, discontinuous (with jumps) L\'evy motion,
are described by two similar differential-integral equations but with different exterior conditions. The non-Gaussianity of the noise manifests as nonlocality at the level of the mean exit time and escape probability.
We consider a numerical approach
for solving these differential-integral (non-local) equations.
A computational analysis is conducted to investigate
the relative importance of jump measure,
diffusion coefficient and non-Gaussianity in affecting mean exit time and escape probability.

\medskip

Our goal is to make mean exit time and escape probability as efficient computational tools, to the applied probability community, for quantifying
stochastic dynamics.

\vskip 12pt

This paper is organized as follows. In section 2, we   recall the generators for L\'evy motions. In
section 3, we  consider SDEs driven by a combination of
Brownian motion and a symmetric $\alpha-$stable process.
Numerical approaches and simulation results are presented
in section 4 and 5, respectively. Finally, the results are summarized in section 6.

\section{L\'evy motion}    \label{motion}

A scalar L\'evy motion is characterized by a linear coefficient
$\theta$, a diffusion parameter $d>0$ and   a non-negative Borel
measure $\nu$, defined on $(\R, \mathcal{B}(\R))$ and concentrated
on $\R \setminus\{0\}$, which satisfies
\begin{equation} \label{levycondition}
  \int_{\R \setminus\{0\} } (y^2 \wedge 1) \; \nu({\rm d}y) < \infty,
\end{equation}
or equivalently
\begin{equation}
  \int_{\R \setminus\{0\} } \frac{y^2}{1+y^2}\; \nu({\rm d}y) < \infty.
\end{equation}
This measure $\nu$ is the so called  the L\'evy jump measure of
the L\'evy motion $L(t)$. We also call $(\theta, d, \nu)$ the
\emph{generating triplet}.

Let $L_t$ be a L\'evy process with the generating triplet
$(\theta, d, \nu)$.
It is known that a scalar L\'evy
motion is completely determined by the L\'evy-Khintchine formula
(See \cite{Apple, Sato-99, PZ}). This says that for any
one-dimensional L\'evy process $L_{t}$, there exists a $\theta \in
R$, $d>0$ and a measure $\nu$ such that
\begin{equation}
Ee^{i \lambda L_{t}}=\exp \{i \theta\lambda t - d t
\frac{\lambda^{2}}{2} +t \int_{\R \setminus \{0\} } (e^{ i \lambda
y}-1 -i \lambda y I_{\{|y| <1 \}} ) \nu({\rm d}y)\},
\end{equation}
where $I_S$ is the indicator function of the set $S$, i.e., it
takes value $1$ on this set and takes zero value otherwise:
$$
I_S(y) =
    \begin{cases}
        1,   &\text{if $y \in S $;}\\
        0,    &\text{if $y \notin S$.}
    \end{cases}
$$

The generator $A$ of the process $L_t$ is defined as
$A  \phi = \lim_{t
\downarrow 0} \frac{P_{t} \phi -\phi}{t}$ where $P_{t} \phi(x)=
E_{x} \phi(L_{t})$ and $\phi$ is any function belonging to the
domain of the operator $A$. Recall that the space $C^2_b( \mathbb{R})$ of $C^2$ functions with bounded derivatives up to order 2 is contained in the domain of $A$, and that for every $\phi\in C^2_b( \mathbb{R})$
(See \cite{Apple, PZ})
\begin{equation} \label{A}
A  \phi(x) =  \theta \phi'(x) + \frac{d}{2} \phi''(x) + \int_{\R
\setminus\{0\}} [\phi(x+ y)-\phi(x) -  I_{\{|y|<1\}} \; y \phi'(x)
] \; \nu({\rm d}y).
\end{equation}

Moreover, the generator for the process   $X_{t} $  in \eqref{sde} is then
\begin{eqnarray} \label{AA0}
A  \phi &=& f(x) \phi^{'}(x)+  \theta \phi'(x) +\frac{d}{2}
   \phi''(x)  \nonumber \\
& &+ \int_{\R \setminus\{0\}} [\phi(x+  y)-\phi(x) -
I_{\{|y|<1\}} \; y \phi'(x) ] \; \nu({\rm d}y).
\end{eqnarray}

For $\alpha\in (0, 2]$, a symmetric $\alpha$-stable
process is a Levy process $L_t$ such that
$$
Ee^{i\lambda L_t}=e^{-t|\lambda|^\alpha}, \quad t>0, \lambda \in \R.
$$
A symmetric 2-stable process is simply a Brownian motion.
When $\alpha\in (0, 2)$,
the generating triplet of the symmetric $\a$-stable process $L_t$
is $(0, 0, \nu_\alpha)$, where
$$
\nu_\a({\rm d}x)=C_\alpha|x|^{-(1+\alpha)}\, {\rm d}x
$$
with $C_\alpha$ given by the formula
$\displaystyle{C_{\alpha} =
\frac{\alpha}{2^{1-\alpha}\sqrt{\pi}}
\frac{\Gamma(\frac{1+\alpha}{2})}{\Gamma(1-\frac{\alpha}{2})}}$.
For more information see \cite{wu, Chen, Schertzer}.

\section{Mean exit time and escape probability}  \label{exit1}

Now we consider the SDE \eqref{sde} with a L\'evy motion
$L_t^{\alpha}$ that has the generating triplet $(0, d, \eps \nu_{\a})$, i.e.,
zero linear coefficient, diffusion coefficient $d\ge 0$ and L\'evy
measure $\eps\nu_\a({\rm d}u)$, with $0<\alpha <2$.
This L\'evy motion is the independent sum of a Brownian motion and a
symmetric $\a$-stable process.
Here $\eps$ is
a non-negative parameter and it does \emph{not} have to be sufficiently small.
Strictly speaking, this is not a $\alpha-$stable L\'evy motion (because the diffusion $d$ may be nonzero), but a L\'evy motion whose jump measure is the same as that of a $\alpha-$stable L\'evy motion.

 We first consider   the mean exit time,  $u(x) \geq 0$,
for an orbit starting at $x$, from a bounded interval $D$.
By the Dynkin formula \cite{Apple, Sato-99} for Markov processes, as in \cite{Naeh, Schuss, Oksendal3}, we  obtain that $u(x) $
satisfies the following differential-integral equation:
\begin{eqnarray}
 A u(x) = -1, \;\;  x\in D \label{exit1D}\\
 u =0, \;\;  x \in D^c, \nonumber
\end{eqnarray}
where the generator $A$ is
\begin{eqnarray} \label{AA}
A  u &=& f(x) u^{'}(x) +  \frac{d}{2}
   u''(x)  \nonumber \\
& &+ \eps \int_{\R \setminus\{0\}} [u(x+  y)-u(x) -
   I_{\{|y|<1\}} \, y u'(x) ] \; \nu_{\a}({\rm d}y),
\end{eqnarray}
and $D^c=\R \setminus D$ is the complement set of $D$.




\bigskip

We also consider the escape probability of a particle whose motion is described by the SDE \eqref{sde}.
The  likelihood of a particle, starting at a point $x$,
first escapes a domain $D$ and   lands
in a subset $E$ of $D^c$ (the complement of $D$) is called escape probability.
This escape probability, denoted by $P_E(x)$, satisfies \cite{liao89,song93} the following equation
\begin{eqnarray}
 A\, P_E(x)&=&0, \quad x \in D, \label{eq.ep}\\
 P_E|_{x \in E}&=&1, \quad  P_E|_{x \in D^c\setminus E}=0, \nonumber
\end{eqnarray}
where $A$ is the generator defined in \eqref{AA}.

In the present paper, we only consider scalar SDEs. For SDEs in higher dimensions, both mean exit time and escape probability will satisfy partial differential-integral equations, and our approaches generally apply.

\section{Numerical schemes}

Noting the principal value of
the integral  $\displaystyle{\int_\R \frac{I_{\{ |y|<\delta\}}(y)
\, y} {|y|^{1+\alpha}}\; {\rm d}y}$  always vanishes for any $\delta>0$,
we will choose the value of $\delta$ in Eq.~(\ref{exit1D})
differently according to the value of $x$.
Eq.~(\ref{exit1D}) becomes
\begin{equation} \label{exit1Dn2}
 \frac{d}{2} \; u''(x)+  f(x)\; u'(x)
  + \eps C_\a  \int_{\R \setminus\{0\}} \frac{u(x+y)-u(x)
     -  I_{\{|y|< \delta \}}(y) \; y u'(x)}{|y|^{1+\alpha}}\; {\rm d}y  = -1,
\end{equation}
for $x\in (a, b)$; and $u(x)=0$ for $x \notin  (a, b)$.

Numerical approaches for the mean exit time and escape probability in the SDEs with Brownian motions were considered in \cite{BrannanDuanErvin, BrannanDuanErvin2}, among others. In the following, we describe the numerical algorithms for the special
case of $(a,b)=(-1,1)$ for clarity of the presentation. The corresponding
schemes for the general case can be extended easily.
Because $u$ vanishes outside $(-1,1)$, Eq.~(\ref{exit1Dn2}) can be
simplified by writing $\int_{\R}=\int_{-\infty}^{-1-x} + \int_{-1-x}^{1-x}
+ \int_{1-x}^{\infty}$,
\begin{eqnarray}
  \frac{d}{2} u''(x) + f(x) u'(x)
  - \frac{\eps C_\a}{\a} \left[\frac{1}{(1+x)^\a}+\frac{1}{(1-x)^\a}\right] u(x) &
\nonumber \\
+ \eps C_\a \int_{-1-x}^{1-x} \frac{u(x+y) - u(x) -
   I_{\{|y|<\delta\}} y u'(x)}{|y|^{1+\a}}\; {\rm d}y & = -1,
\label{exit1Dn3}
\end{eqnarray}
for $x \in (-1,1)$; and $u(x)=0$ for $x \notin  (-1, 1)$.

   Noting $u$ is not smooth at the boundary points $x=-1, 1$,
in order to ensure the integrand is smooth, we rewrite Eq.~(\ref{exit1Dn3})
as
\begin{eqnarray}
  \frac{d}{2} u''(x) + f(x) u'(x)
  - \frac{\eps C_\a}{\a} \left[\frac{1}{(1+x)^\a}+\frac{1}{(1-x)^\a}\right] u(x) &
\label{exit1Dn4}\\
+ \eps C_\a \int_{-1-x}^{-1+x} \frac{u(x+y) - u(x)}{|y|^{1+\a}}\; {\rm d}y
  + \eps C_\a \int_{-1+x}^{1-x} \frac{u(x+y)-u(x)-y u'(x)}{|y|^{1+\a}}\; {\rm d}y & = -1,
\nonumber
\end{eqnarray}
for $x\geq 0$, and
\begin{eqnarray}
  \frac{d}{2} u''(x) + f(x) u'(x)
  - \frac{\eps C_\a}{\a} \left[\frac{1}{(1+x)^\a}+\frac{1}{(1-x)^\a}\right] u(x) &
\label{exit1Dn5}\\
+ \eps C_\a \int_{1+x}^{1-x} \frac{u(x+y) - u(x)}{|y|^{1+\a}}\; {\rm d}y
  + \eps C_\a \int_{-1-x}^{1+x} \frac{u(x+y)-u(x)-y u'(x)}{|y|^{1+\a}}\; {\rm d}y & = -1,
\nonumber
\end{eqnarray}
for $x < 0$. We have chosen $\delta= \text{min}\{ |-1-x|,|1-x|\}$.

   Let's divide the interval $[-2,2]$ into $4J$ sub-intervals
and define $x_j=jh$ for $-2J\leq j \leq 2J$ integer, where $h=1/J$. We denote
the numerical solution of $u$ at $x_j$ by $U_j$. Let's discretize
the integral-differential equation (\ref{exit1Dn4}) using
central difference for derivatives and ``punched-hole'' trapezoidal
rule
\begin{equation}
  \begin{split}
  & \frac{d}{2} \frac{U_{j-1} - 2U_j + U_{j+1}}{h^2}
   + f(x_j) \frac{U_{j+1} - U_{j-1}}{2h}
   -  \frac{\eps C_\a}{\a} \left[\frac{1}{(1+x_j)^\a}+\frac{1}{(1-x_j)^\a}\right] U_j \\
  & + \eps C_\a h \sum^{-J+j}_{k=-J-j}\!\!\!\!\!\!{''} \;
    {\frac{U_{j+k} - U_j}{|x_k|^{1+\a}} }
   + \eps C_\a h \sum^{J-j}_{k=-J+j,k\neq 0}\!\!\!\!\!\!\!\!\!{''} \;
    {\frac{U_{j+k} - U_j -(U_{j+1}-U_{j-1}) x_k/2h}{|x_k|^{1+\alpha}} } = -1,
  \end{split}
 \label{nm1D1a}
\end{equation}
where $j = 0,1,2, \cdots, J-1$. The modified summation symbol $\sum{''}$
means that the quantities corresponding to the two end summation
indices are multiplied by $1/2$.
\begin{equation}
  \begin{split}
  & \frac{d}{2} \frac{U_{j-1} - 2U_j + U_{j+1}}{h^2}
    + f(x_j) \frac{U_{j+1} - U_{j-1}}{2h}
   -  \frac{\eps C_\a}{\a} \left[\frac{1}{(1+x_j)^\a}+\frac{1}{(1-x_j)^\a}\right] U_j \\
  & + \eps C_\a h \sum^{J-j}_{k=J+j}\!\!\!\!{''} \;
    {\frac{U_{j+k} - U_j}{|x_k|^{1+\a}} }
   + \eps C_\a h \sum^{J+j}_{k=-J-j,k\neq 0}\!\!\!\!\!\!\!\!\!{''} \;
    {\frac{U_{j+k} - U_j -(U_{j+1}-U_{j-1}) x_k/2h}{|x_k|^{1+\alpha}} } = -1,
  \end{split}
 \label{nm1D1b}
\end{equation}
where $j = -J+1, \cdots, -2,-1$.
The boundary conditions require that
the values of $U_j$ vanish if the index $|j|\geq J$.

The truncation errors of the central differencing schemes for derivatives in
(\ref{nm1D1a})  and (\ref{nm1D1b}) are of 2nd-order $O(h^2)$.
From the error analysis
of Navot (1961) \cite{Navot61}, the leading-order error of the
quadrature rule is $-\zeta(\alpha-1) u''(x) h^{2-\alpha} +
O(h^2)$, where $\zeta$ is the Riemann zeta function. Thus,
the following scheme have 2nd-order accuracy
for any $0<\a<2$,
 $j = 0,1,2, \cdots, J-1$
\begin{equation}
  \begin{split}
  & C_h \frac{U_{j-1} - 2U_j + U_{j+1}}{h^2}
   + f(x_j) \frac{U_{j+1} - U_{j-1}}{2h}
   -  \frac{\eps C_\a}{\a} \left[\frac{1}{(1+x_j)^\a}+\frac{1}{(1-x_j)^\a}\right] U_j \\
  & + \eps C_\a h \sum^{-J+j}_{k=-J-j}\!\!\!\!\!\!{''} \;
    {\frac{U_{j+k} - U_j}{|x_k|^{1+\a}} }
   + \eps C_\a h \sum^{J-j}_{k=-J+j,k\neq 0}\!\!\!\!\!\!\!\!\!{''} \;
    {\frac{U_{j+k} - U_j -(U_{j+1}-U_{j-1}) x_k/2h}{|x_k|^{1+\alpha}} } = -1,
  \end{split}
 \label{nm1D2a}
\end{equation}
where $\displaystyle{C_h = \frac{d}{2} -
\eps C_\a \zeta(\alpha-1) h^{2-\a}}$.
Similarly,
for $j = -J+1, \cdots, -2,-1$,
\begin{equation}
  \begin{split}
  & C_h \frac{U_{j-1} - 2U_j + U_{j+1}}{2h^2}
    + f(x_j) \frac{U_{j+1} - U_{j-1}}{2h}
   -  \frac{\eps C_\a}{\a} \left[\frac{1}{(1+x_j)^\a}+\frac{1}{(1-x_j)^\a}\right] U_j \\
  & + \eps C_\a h \sum^{J-j}_{k=J+j}\!\!\!\!{''} \;
    {\frac{U_{j+k} - U_j}{|x_k|^{1+\a}} }
   + \eps C_\a h \sum^{J+j}_{k=-J-j,k\neq 0}\!\!\!\!\!\!\!\!\!{''} \;
    {\frac{U_{j+k} - U_j -(U_{j+1}-U_{j-1}) x_k/2h}{|x_k|^{1+\alpha}} } = -1,
  \end{split}
 \label{nm1D2b}
\end{equation}
where $j = -J+1, \cdots, -2,-1,0,1,2, \cdots, J-1$.
$U_j=0$ if $|j|\geq J$.

We solve the linear system (\ref{nm1D2a}-\ref{nm1D2b})
by direct LU factorization or
the Krylov subspace iterative method GMRES.

We find that the desingularizing term ($I_{\{|y|<\delta\}} y u'(x)$)
does not have any effect on the numerical results, regardless whether
we use LU or GMRES for solving the linear system. In this case,
we can discretize the following equation instead of \eqref{exit1Dn3}
\begin{eqnarray}
  \frac{d}{2} u''(x) + f(x) u'(x)
  - \frac{\eps C_\a}{\a} \left[\frac{1}{(1+x)^\a}+\frac{1}{(1-x)^\a}\right] u(x) &
\nonumber \\
+ \eps C_\a \int_{-1-x}^{1-x} \frac{u(x+y) - u(x)}
     {|y|^{1+\a}}\; {\rm d}y & = -1,
\label{exit1Dn6}
\end{eqnarray}
where the integral in the equation is taken as Cauchy principal value integral.
Consequently, instead of \eqref{nm1D2a} and \eqref{nm1D2b},
we have only one discretized equation for any $0<\a<2$ and
$j = -J+1, \cdots, -2,-1,0,1,2, \cdots, J-1$
\begin{equation}
  \begin{split}
  & C_h \frac{U_{j-1} - 2U_j + U_{j+1}}{h^2}
   + f(x_j) \frac{U_{j+1} - U_{j-1}}{2h} \\
  & -  \frac{\eps C_\a U_j}{\a} \left[\frac{1}{(1+x_j)^\a}+\frac{1}{(1-x_j)^\a}\right]
   + \eps C_\a h \sum^{J-j}_{k=-J-j,k\neq 0}\!\!\!\!\!\!\!\!\!{''} \;
    {\frac{U_{j+k} - U_j}{|x_k|^{1+\alpha}} } = -1.
  \end{split}
 \label{nm1D3}
\end{equation}

\bigskip

With minor changes, we also have the  scheme for simulating escape probability characterized by the equation
\eqref{eq.ep}.

\section{Numerical results}

\subsection{Verification}
\subsubsection{Comparing with analytical solutions}

  In order to verify that the numerical integration scheme
for treating the improper integral in \eqref{exit1D}
is implemented correctly,
we compute the left-hand side(LHS) of \eqref{exit1D}
by substituting $u(x)=1-x^2$, $d=0$, $f(x)\equiv 0$
$\eps=1$ and $(a,b)=(-1,1)$
\begin{equation}
  LHS = \begin{cases}
 -C_{\alpha}\left[(1+x)^{1-\a} \left(
   \frac{1-x}{\a} - \frac{2x}{1-\a} + \frac{1+x}{2-\a}\right)
      \right.  & \\
   \left.  + (1-x)^{1-\a} \left( \frac{1+x}{\a}
    + \frac{2x}{1-\a} + \frac{1-x}{2-\a}\right) \right]
    &  \text{if $\alpha\neq 1$;} \\
 -C_{\alpha}\left(4 + 2x\ln\frac{1-x}{1+x}\right),
    &  \text{if $\alpha = 1$.}
        \end{cases}
\label{eq.lhs}
\end{equation}

\begin{figure}
\begin{center}
\vspace*{-1.1in}
\includegraphics*[width=\linewidth]{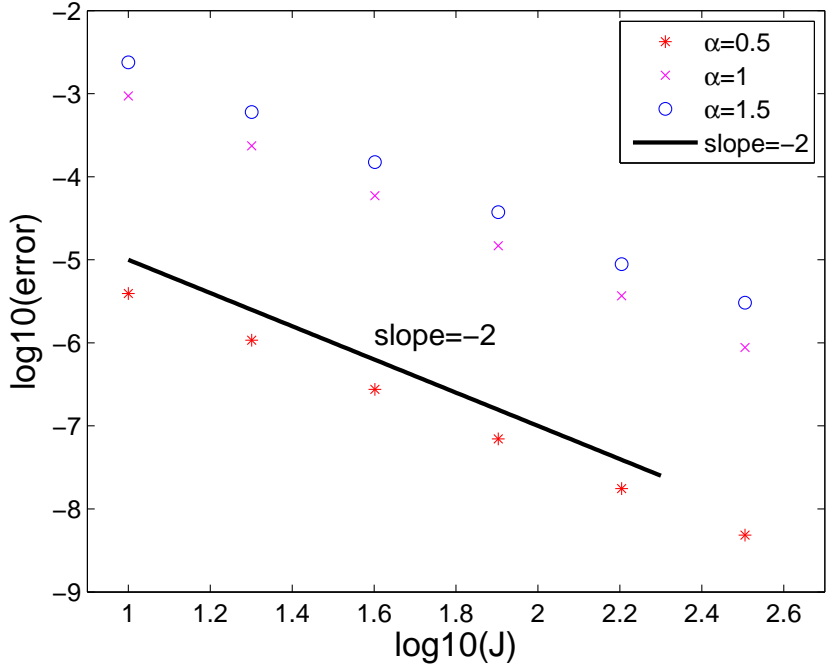}
\end{center}
\caption{The error of the numerical values of the left-hand side
of Eq.~\eqref{exit1D} compared with the analytical expression
in \eqref{eq.lhs} for $u(x)=1-x^2$, $d=0$, $f(x)\equiv 0$, $\eps=1$
and $(a,b)=(-1,1)$.
The results are computed at $x=-0.5$ for different values
of $\a=0.5$ (marked by $*$'s), $1$
(the x's) and $1.5$ (the o's) and different resolutions $J=10, 20, 40,
80, 160$ and $320$. Also shown is an illustrating solid line with slope equal
to $-2$.
}
\label{lhs}
\end{figure}
   Figure~\ref{lhs} shows the differences (the errors) between
the numerical and the analytical values of LHS of \eqref{exit1D}
at the fixed value $x=-0.5$ for different resolutions $J=10, 20, 40,
80, 160$ and $320$. We plot $\log_{10}(error)$ against $\log_{10}(J)$
for $\a=0.5, 1, 1.5$,
where $h=1/J$ and $error$ is the difference between the numerical
and the analytical values of LHS. Clearly, the numerical results
show that the error of computing LHS decays as $O(h^2)$. The second-order
accuracy is expected from the error analysis of the numerical integration
method \eqref{nm1D3}. For a fixed resolution $h$, the error increases
as $\alpha$ increases due to the fact that LHS in \eqref{eq.lhs}
becomes more singular at $x=-1$ and $1$ as $\alpha$ increases.

\begin{figure}
\begin{center}
\vspace*{-0.8in}
\includegraphics*[width=\linewidth]{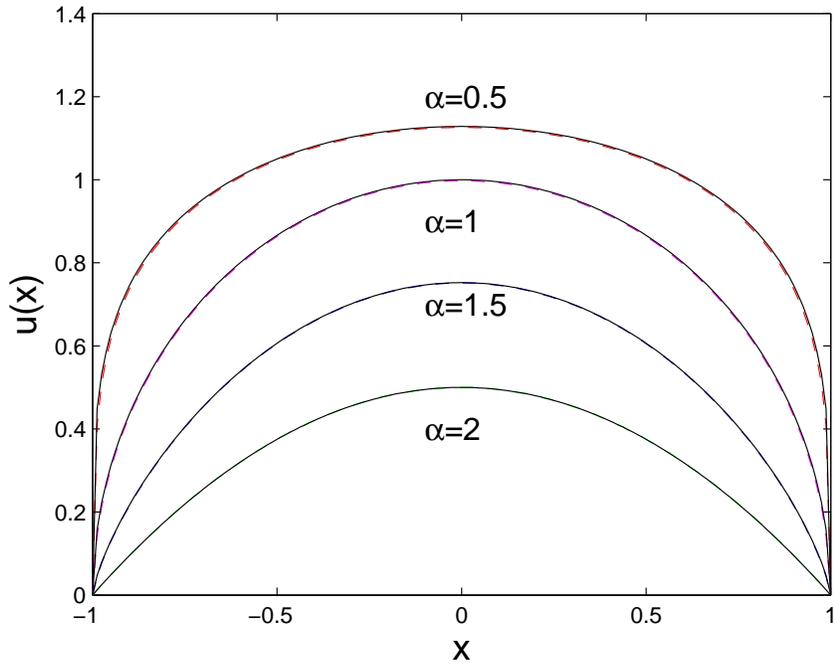}
\end{center}
\caption{The mean exit time $u(x)$ for the special case
of the symmetric $\a$-stable process
$d=0$, $f(x)\equiv 0$, $\eps=1$ and $(a,b)=(-1,1)$.
The dashed lines are the numerical solutions
for different values of $\a=0.5, 1, 1.5$ obtained by solving
the discretized equations \eqref{nm1D3} with
the resolution $J=80$, while the solid lines represent the corresponding
analytical solutions \eqref{eq.as} including $\a=2$.}
\label{comp_sol}
\end{figure}
Next, we compare the numerical solution with the analytical solution
for the mean exit time \cite{Getoor}
\begin{equation}
u(x) =
\frac{\sqrt{\pi}(b^2-x^2)^{\alpha/2}}
{2^{\alpha}\Gamma(1+\alpha/2)\Gamma(1/2+\alpha/2)}
\label{eq.as}
\end{equation}
in the special case of \eqref{exit1D}
in which $d=0$, $f(x)\equiv 0$, $\eps=1$ and $(a,b)=(-b,b)$ with $b>0$.
Figure~\ref{comp_sol} shows the numerical solutions (the dashed lines)
obtained by solving the discretized equations \eqref{nm1D2a} and \eqref{nm1D2b}
or \eqref{nm1D3} with the fixed resolution $J=80$, $(a,b)=(-1,1)$ and different
values of $\a=0.5, 1, 1.5$, while the corresponding analytical solutions
are shown with the solid lines. The comparison shows that
the numerical solutions are very accurate
as one can hardly distinguish the numerical solution from the corresponding
analytical one. Note that, in this case of $(a,b)=(-1,1)$,
for a fixed value of the starting point $x$,
the mean exit time $u(x)$ decreases when $\alpha$ increases in the interval
$(0,2]$. Later, we will see the dependence of the mean exit time on $\alpha$
is much more complicated when the size of the interval $b-a$ is increased.

\begin{figure}
\begin{center}
\includegraphics*[width=\linewidth]{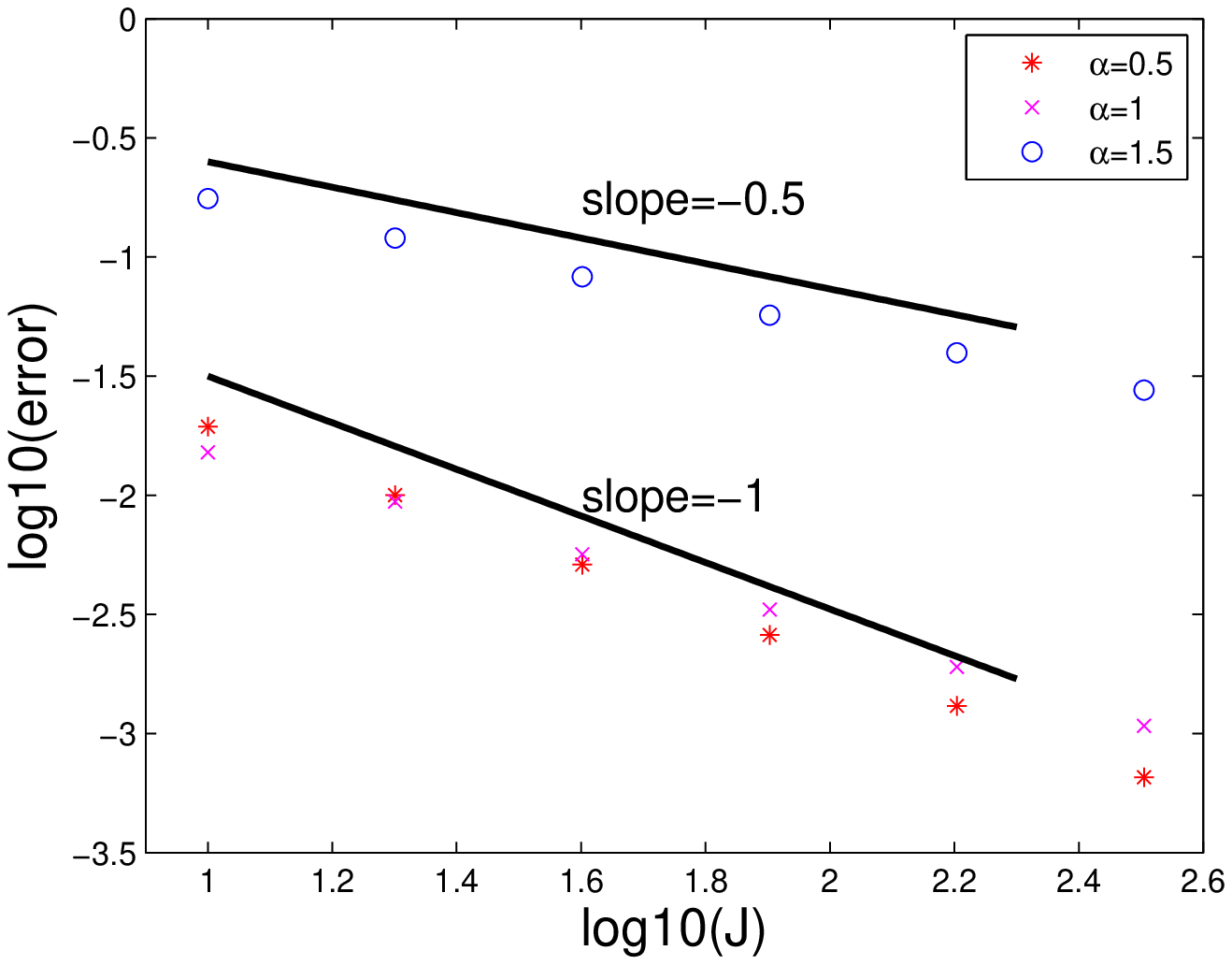}
\end{center}
\caption{The error of the numerical solution of the mean-exit time $u(x)$
from the discretized equations \eqref{nm1D1a} and \eqref{nm1D1b}
using the "punched-hole" trapezoidal rule
for $d=0$, $f(x)\equiv 0$, $\eps=1$ and $(a,b)=(-1,1)$.
The results are computed at $x=-0.5$ for different
values of $\a=0.5$ (marked by $*$'s), $1$ (the x's) and $1.5$ (the
o's) and different resolutions $J=10, 20, 40, 80, 160$ and $320$.
Also shown are two illustrating solid lines with slope equal to $-1$
and $-0.5$ respectively. }
\label{aswoc}
\end{figure}
Figure~\ref{aswoc} gives the error in the numerical solution
to the discretized equations \eqref{nm1D1a} and \eqref{nm1D1b} derived
by using "punched-hole" trapezoidal rule.
Here, we compute $error= |u(-0.5)-U(-0.5)|$ by comparing
the mean exit time at $x=-0.5$
for different resolutions and values of $\a$,
where $u$ and $U$ denotes the analytical
and the numerical solution respectively.
For a fixed resolution, the numerical error has similar sizes
for $\a=0.5$ and $\a=1$ but it is much larger in the case of $\a=1.5$.
The analysis shows that the rate of decay in the error
as the resolution increases is $O(h^{2-\a})$.
Our numerical results in Fig.~\ref{aswoc}
show slower decaying rates than those in the theory for $\a=0.5$ and $1$,
which is due to the non-smoothness of the solution at $x=-1,1$.

\begin{figure}
\begin{center}
\vspace*{-1.1in}
\includegraphics*[width=\linewidth]{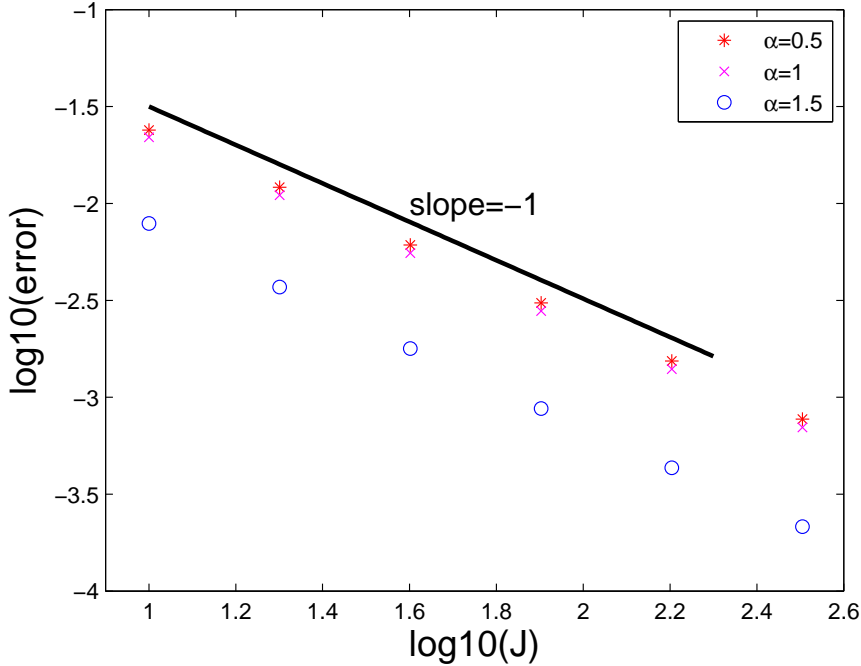}
\end{center}
\caption{The error of the numerical solution of the mean-exit time $u(x)$
from the discretized equations \eqref{nm1D2a} and \eqref{nm1D2b}
or \eqref{nm1D3} with the correction term
for $d=0$, $f(x)\equiv 0$, $\eps=1$ and $(a,b)=(-1,1)$.
The results are computed at $x=-0.5$ for different
values of $\a=0.5$ (marked by $*$'s), $1$ (the x's) and $1.5$ (the
o's) and different resolutions $J=10, 20, 40, 80, 160$ and $320$.
Also shown is an illustrating solid line with slope equal to $-1$. }
\label{aswc}
\end{figure}

Figure~\ref{aswc} is the same as Fig.~\ref{aswoc} except that the numerical
results are obtained from
the discretized equations \eqref{nm1D2a} and \eqref{nm1D2b}
or \eqref{nm1D3} with the correction term that removes the leading-order
quadrature error. Although the numerical analysis predicts the decaying
rate of the numerical error is $O(h^2)$, the numerical
results shown in Fig.~\ref{aswc} indicate the rate of decay is only $O(h)$,
because the analytical solution $u(x)$ given in \eqref{eq.as} has infinite
derivatives at $x=-1$ and $1$. Note that we have demonstrated
in Fig.~\ref{lhs} that the convergence order would be 2 if the solution
$u$ were smooth on the whole closed interval $[-1,1]$.
Though the convergence order is only 1, it become independent of $\a$
after we add the correction term and the numerical error is two orders
of magnitude smaller than that without the correction term when $\alpha=1.5$.

\subsubsection{Comparing with asymptotic solutions}

 Since we do not have analytical solutions in closed form
for the mean first exit time $u$
to Eq.~\eqref{exit1D} in the general case of $d\neq 0, f\neq 0$,
we calculate the asymptotic solutions
for small values of $\eps$, i.e., small pure jump measure
in the L\'evy motion. Then, we test our numerical schemes by comparing
the results with those from the corresponding asymptotic solutions.

We look for solution to Eq.~\eqref{exit1D}
in the form of expansion
\begin{equation}
u(x)=u_0(x)+\varepsilon u_1(x)+ \varepsilon^2
u_2(x)+\cdots.
\label{eq.uexp}
\end{equation}
By setting $\varepsilon=0$ in Eq.~\eqref{exit1D},
we obtain the equation for $u_0$
\begin{equation}
\frac{d}{2}u_0''(x) + f(x) u_0'(x)=-1.
\end{equation}
The solution is given by
\begin{equation}
u_0(x)=\int^x v(y)\, {\rm d}y + A_2, \quad \text{where }
v(x) = \frac{-\frac{2}{d} \int^x
e^{\frac{2}{d} \int^z f(y)\,{\rm d}y}\,{\rm d}z +A_1}
{ e^{\frac{2}{d}\int^x f(y)\,{\rm d}y}},
\label{eq.as0}
\end{equation}
where $A_1$ and $A_2$ are integration constants that can be determined
by the boundary conditions.
Substituting the expansion \eqref{eq.uexp} into \eqref{exit1D}
and discarding the terms with second or higher powers of $\eps$,
we obtain the equation for $u_1$
\begin{equation}
\frac{d}{2}u_1''(x) + f(x) u'_1(x)
+ C_{\alpha} \int_{\mathbb{R}\setminus \{0\}}
\frac{u_0(x+y)-u_0(x)- I_{\{|y|<1\}}(y)
yu_0'(x)}{|y|^{1+\alpha}}\,{\rm d}y=0.
\end{equation}
Denoting $\displaystyle{g(x):=C_{\alpha}
\int_{\mathbb{R}\setminus \{0\}}
\frac{u_0(x+y)-u_0(x)- I_{\{|y|<1\}}(y)
yu_0'(x)}{|y|^{1+\alpha}}\,{\rm d}y}$, we have
\begin{equation}
u_1(x)=\int^x w(y)\,{\rm d}y + A_4, \quad \text{where }
w(x)=\frac{-\frac{2}{d}\int^x e^{\int^z \frac{2}{d} f(y)
\,{\rm d}y} g(z) \,{\rm d}z + A_3}{e^{\frac{2}{d} \int^x f(y)\,{\rm d}y}},
\label{eq.as1}
\end{equation}
where $A_3$ and $A_4$ are integration constants that can be determined
by the boundary conditions.

Let us consider the special case $f(x)\equiv 0$, $d=1$ and $(a,b)=(-1,1)$.
According to the general solution \eqref{eq.as0} and \eqref{eq.as1},
the zeroth and first-order solutions are
\begin{eqnarray}
   u_0(x)&=&1-x^2, \label{eq.as2}\\
u_1(x)&=& \begin{cases} \frac{2C_{\alpha}}{\a(1-\a)(2-\a)(3-\a)(4-\a)}
    \left[(2-\a) (1+x)^{4-\a} \right. & \nonumber \\
   + (4-\a) (1+x)^{3-\a} (1-x) + (4-\a) (1+x) (1-x)^{3-\a} & \\
   \left. + (2-\a) (1-x)^{4-\a}
    - 2^{4-\a} (2-\a) \right], & \text{for }\a \in (0,1) \cup (1,2),\\
  \frac{2C_{\alpha}}{3}\left[2x^2 -2 - 4\ln2
       + (2+x) (1-x)^2 \ln(1-x)\right. & \\
   \left. + (2-x) (1+x)^2 \ln(1+x) \right], & \text{for }\a=1.
    \end{cases}
\end{eqnarray}

\begin{figure}
\begin{center}
\includegraphics*[width=\linewidth]{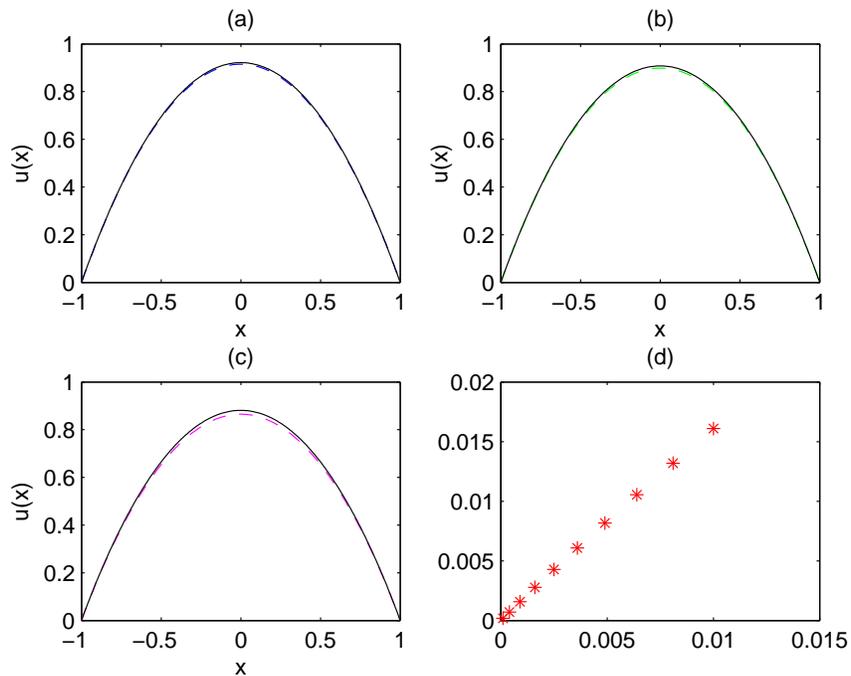}
\end{center}
\caption{Comparison between the numerical solution
to \eqref{exit1D} and the asymptotic solution \eqref{eq.as2}
for small $\eps$ with $f\equiv 0$, $d=1$ and $(a,b)=(-1,1)$.
(a) $\eps=0.1$ and $\a=0.5$. The numerical solution
is displayed with dashed line
while the asymptotic solution $u_0+\varepsilon u_1$
is shown with solid black line. (b) Same as (a) except $\a=1$.
(c) Same as (a) except $\a=1.5$. (d) The difference between
the numerical solution $U(0)$ and the asymptotic solution
$u_0(0)+\varepsilon u_1(0)$ is plotted against $\eps^2$
for $\a=1.5$ and the fixed resolution $J=100$.}
\label{comp_asy}
\end{figure}
To further verify our numerical methods, we compare the numerical solution
to \eqref{exit1D} with the asymptotic solution \eqref{eq.as2}
$u_0+\varepsilon u_1$ for the case $\eps=0.1$, $f\equiv 0$, $d=1$ and
$(a,b)=(-1,1)$, as shown in Fig.~\ref{comp_asy}.
The plots show that the two solutions are very close for $\a=0.5$
(Fig.~\ref{comp_asy}(a)), $\a=1$ (Fig.~\ref{comp_asy}(b))
and $\a=1.5$ (Fig.~\ref{comp_asy}(c)). Furthermore, Fig.~\ref{comp_asy}(d)
shows that the difference between the two solutions is proportional
to $\eps^2$, which is expected as the asymptotic solution
$u_0 +\eps u_1$ given by \eqref{eq.as2} is only accurate up to $O(\eps)$.

\subsection{Dependence of the mean exit time on the size of domain}

\subsubsection{Pure jump: $d=0, f\equiv 0$}
\begin{figure}
\begin{center}
\vspace*{-1.0in}
\includegraphics*[width=13.5cm]{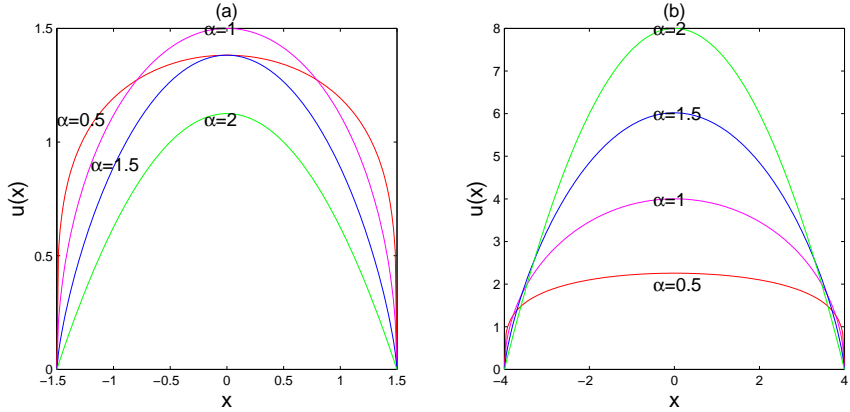}
\end{center}
\caption{Dependence of the mean exit time $u(x)$
on the size of domain $D$ for pure $\a$-stable L\'evy motion,
i.e., $d=0$, $f(x)\equiv 0$, $\eps=1$. (a) The mean first exit time $u(x)$
from the analytical solution \eqref{eq.as}
for the domain $D=(-1.5,1.5)$ and $\a=0.5, 1, 1.5, 2$.
(b) Same as (a) except $D=(-4,4)$.
}
\label{dds}
\end{figure}

It is well-known that the $\a$-stable L\'evy process has larger jumps
with lower jump frequencies for small values of $\a$ ($0<\a<1$)
while it has smaller jumps with higher jump probabilities for values
of $\a$ closer to 2. Given a bounded domain $D$, it would be interesting
to know, for symmetric $\a$-stable L\'evy motion ($d=0$ and $f\equiv 0$),
the exit times out of domain $D$ are shorter for small values
of $\a$ or large values of $\a$. The answer, given
by the analytic solution \eqref{eq.as},
depends on the size of domain $D$. To illustrate the dependence,
we compare the solutions for four representative values
of $\a$, i.e., $0.5, 1, 1.5$, and $2$.
For a small-size domain $D=(-b,b)$ with
$0<b\leq 1.25$, it is easier to leave the domain for larger values of $\a$
and any starting point $x$, such as the case $b=1$
shown in Fig.~\ref{comp_sol}. In contrast,
for large domains $D$ with $b\geq 3$, the exit times are
shorter when $\a$ is smaller except for starting points near the boundary
such as the case $b=4$ shown in Fig.~\ref{dds}(b).
For median-sized domains $D=(-b,b)$ with $1.25<b<3$,
whether the mean exit time is shorter for smaller or larger values of $\a$
depends on both the position of the starting point $x$
and the size of symmetric domain $2b$,
such as the case $b=1.5$ shown in Fig.~\ref{dds}(a).

The results on the mean exit time show that one has to consider both the domain
size and the value of $\a$ when deciding which of the two competing
factors in $\a$-stable L\'evy motion, the jump frequency or the jump size,
is dominant. The small jumps with high frequency, corresponding
to L\'evy motion with $\a$ closer to $2$, make it easier
to exit the small domains.
On the other hand, it is easier to exit large domains for
$\a$-stable L\'evy motion with $\a$ closer to $0$,
which has the characteristics of the large jumps with low frequency.
Another observation from the results is that, for smaller values
of $\a$ ($0<\a<1$), the mean exit time profiles becomes flatter
away from the boundary, such as the graphs for $\a=0.5$
in Fig.~\ref{comp_sol} and Fig.~\ref{dds}(b).
This implies that the jump sizes of the processes
are usually larger than the domain sizes, thus the mean exit times
have small variations for different starting positions.

\subsubsection{Ornstein-Uhlenbeck(O-U) potential: $f(x) = -x$}

In the deterministic case ${\rm d} X_t = - X_t \, {\rm d}t$,
the origin is the sole stable point and the particle is driven
toward the origin with the velocity proportional to its distance
to the unique stable point. When a particle is subject to
the L\'evy motion \eqref{sde} defined in the beginning of Sec.~\ref{exit1}
the mean exit time out of a bounded domain becomes finite.
We emphasize that, in this paper, we consider the L\'evy motion
defined by the generator in \eqref{AA} where we vary the diffusion
coefficient $d$ and the parameter $\eps$ independently.
\begin{figure}
\begin{center}
\vspace*{-1.0in}
\includegraphics*[width=14.0cm]{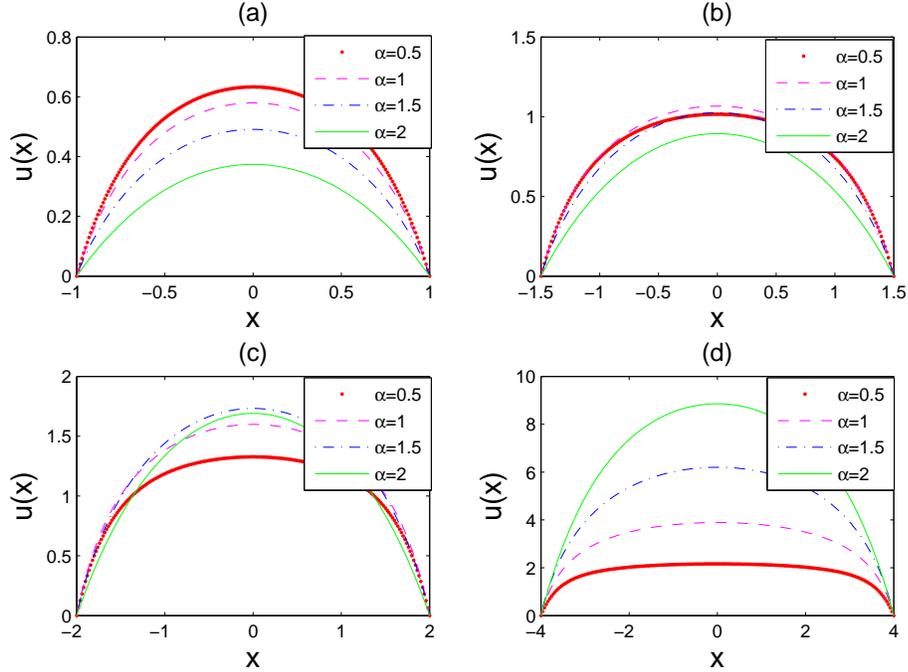}
\end{center}
\caption{Dependence of the mean exit time $u(x)$
on the size of domain $D$ for the case of Ornstein-Uhlenbeck potential
($f(x)=-x$) with both Gaussian ($d=1$) and non-Gaussian ($\eps=1$) noises.
(a) The mean first exit time $u(x)$
for $D=(-1,1)$ and $\a=0.5, 1, 1.5, 2$.
(b) Same as (a) except $D=(-1.5,1.5)$.
(c) Same as (a) except $D=(-2,2)$.
(d) Same as (a) except $D=(-4,4)$.
}
\label{ddsoud1}
\end{figure}
Figure~\ref{ddsoud1} illustrates
the dependence of the exit time $u$
on the size of the domain for the case $f(x)=-x$, $d=1$, $\eps=1$
and the four typical values of $\a=0.5, 1, 1.5, 2$.
Figure~\ref{ddsoud1}(a) shows that, for small domains
$D=(-b,b)$ such as $b=1$ and from any starting point $x$ in $D$,
the mean exit times are shorter for larger values of $\a$,
in agreement with the corresponding result in the absence of
the driving force $f=0$ and the Gaussian noise $d=0$ shown
in Fig.~\ref{comp_sol}.
On the other hand, for large domains such as $D=(-b,b)$ with $b=4$,
Fig.\ref{ddsoud1}(d) demonstrates that the mean exit times
are longer for larger values of $\a$. Again, the behavior agrees in general
with the results in the previous pure jump case with $f=d=0$ and $\eps=1$
shown in Fig.~\ref{dds}(b).
Note that, for the starting points
near the boundaries, the relations between the mean exit times and
the values of $\a$ are different in these two cases:
the mean exit times increases as $\a$ is raised for the case
of nonzero Gaussian noise $d=1$ and the nonzero driving force $f(x)=-x$
while the pure jump L\'evy motion ($f=0$ and $d=0$) has the opposite behavior.
The dependence of the mean exit times on the value of $\a$ is mixed
for median-sized domains, such as $D=(-1.5,1.5)$ and $D=(-2,2)$
shown in Fig.\ref{ddsoud1}(b) and (c) respectively.

\begin{figure}
\begin{center}
\vspace*{-1.0in}
\includegraphics*[width=14.0cm]{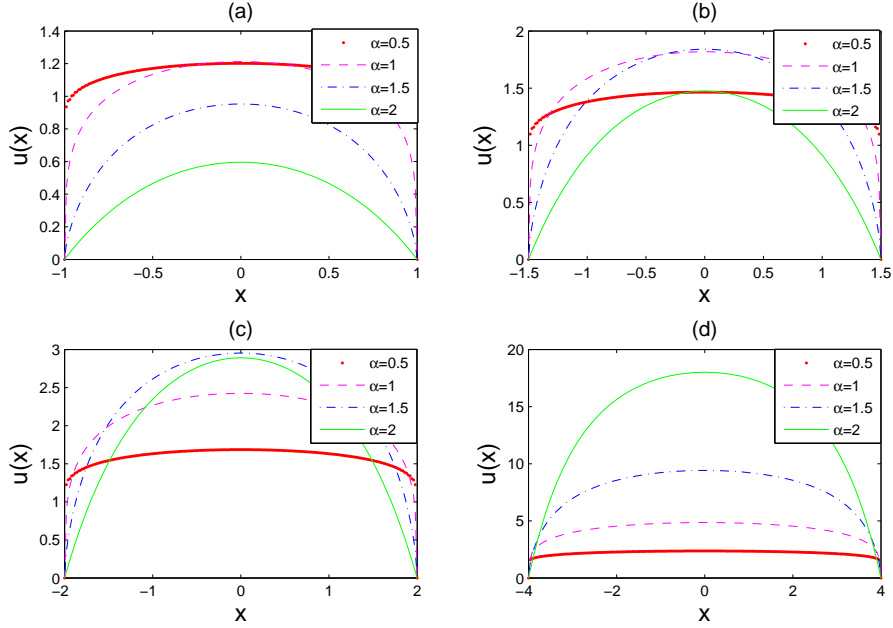}
\end{center}
\caption{Dependence of the mean exit time $u(x)$
on the size of domain $D$ for the case of Ornstein-Uhlenbeck potential
($f(x)=-x$) with pure non-Gaussian noises $d=0$ and $\eps=1$.
(a) The mean first exit time $u(x)$
for $D=(-1,1)$ and $\a=0.5, 1, 1.5, 2$.
(b) Same as (a) except $D=(-1.5,1.5)$.
(c) Same as (a) except $D=(-2,2)$.
(d) Same as (a) except $D=(-4,4)$.
}
\label{ddsoud0}
\end{figure}
Figure~\ref{ddsoud0} shows the mean first exit times when Gaussian noise
is removed while keeping other factors the same, i.e., $f(x)=-x$,
$\eps=1$ but $d=0$. The dependence on the size of the domain is similar to
the previous case with Gaussian noise $d=1$.
However, as shown in Fig.~\ref{ddsoud0}, we find that, in the presence
of the O-U potential ($f(x)=-x$)
and without Gaussian noise ($d=0$), the mean first exit time $u(x)$
not only has a flat profile also is discontinuous at the boundaries $x=\pm b$
for $\a=0.5$. We find that
it is also true for other values of $\a$ in $(0,1)$.  A possible explanation for this is as follows: When $\alpha \in (0, 1)$, the original first order differential operator plays the dominant role, while when $\alpha\in (1, 2)$, the integral operator plays the dominant role.
To obtain the discontinuous numerical solutions in these cases,
we replace the central differencing scheme for the term $f(x) u'(x)$
in \eqref{nm1D3} with a second-order one-sided difference
for the first and last interior grid point.

\subsection{Effect of the noises}

\begin{figure}
\begin{center}
\vspace*{-1.0in}
\includegraphics*[width=14.0cm]{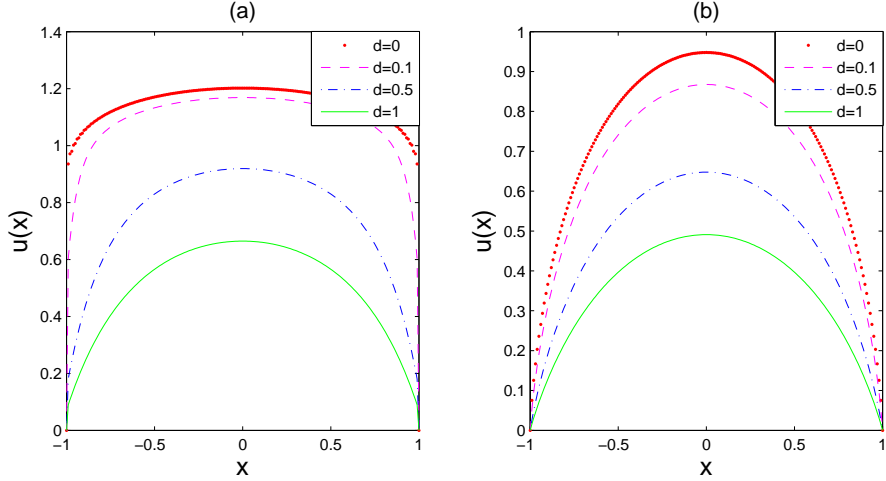}
\end{center}
\caption{Effect of Gaussian noise on the mean exit times
of the L\'evy motion $L_t^\a$. (a) The mean exit times
for $f(x)=-x, \eps=1, \a=0.5, D=(-1,1)$ and $d=0, 0.1, 0.5, 1$.
(b) Same as (a) except $\a=1.5$.
}
\label{egn}
\end{figure}
Having discussed the dependence of the mean exit times on the value of $\a$
in the L\'evy motion $L_t^\a$ defined in the beginning of Sec.~\ref{exit1},
we examine the effect
of Gaussian noise on the profile of the mean exit time $u(x)$ as a function
of the location of the starting point $x$.
We vary the values of the diffusion coefficient $d$ and the parameter
$\eps$ in L\'evy measure independently. Thus, $L_t^\a$ in the SDE,
as defined by the generator in \eqref{AA}, is
different than the traditional L\'evy motion where the diffusion coefficient
$d$ and the coefficient $\eps$ in L\'evy measure are changed in tandem.

Consider the L\'evy motion driven by O-U potential $f(x) = -x$
and the symmetric domain $D=(-1,1)$.
Figure~\ref{egn} shows the numerical results of the mean exit times with
$\eps=1$ for different amount of Gaussian noises $d=0, 0.1, 0.5, 1$.
For small values of $\a$ ($0<\a<1$), such as $\a=0.5$ as shown in Fig.~\ref{egn}(a), the mean exit time shapes as function of the starting point $x$
 change dramatically as the amount of Gaussian noises $d$ increases.
For small amount of Gaussian noises, the mean exit time profile is flat
in the middle of the domain and drops to zero quickly near the boundary points;
for large amount of Gaussian noises, the mean exit time profile become
more parabla-like as shown in the graph for $d=1$.

It is worth pointing out
that the mean first exit time is discontinuous at the boundary $x=\pm 1$
in the case of pure non-Gaussian noise $d=0$ and   $\alpha\in (0, 1)$,
i.e., the limits $\lim_{x\rightarrow \pm 1} u(x)$ are nonzero
while $u(\pm 1)=0$. As mentioned in previous section,
we have to use an one-sided difference scheme near the boundary
to avoid numerically differentiating across discontinuities.
From our numerical simulations for other values of $\a$ and domain sizes
(not shown here), we find that the mean exit time $u(x)$ driven by
O-U potential with "pure" $\a$-stable jump only
and $0<\a<1$ would be discontinuous at the boundary of the domain.
Recall that, in the absence of deterministic
driving force ($f\equiv 0$) and Gaussian noise ($d=0$),
the mean exit time profile given in \eqref{eq.as} become more "discontinuous"
at the boundary as $\a \rightarrow 0+$ (more precisely, the derivative
of $u(x)$ goes to infinity faster near the boundary
for smaller values of $\a$). Our numerical results show that
adding O-U potential would cause the mean exit times
be discontinuous at the boundary of the domain for all values of $\a$
in $(0,1)$ and any domain size.

For large values of $\a$ ($1\leq \a \leq 2$), such as $\a=1.5$
as shown in Fig.~\ref{egn}(a), the mean exit time shapes are similar
to the parabolic shape as in the pure Gaussian noise case. Clearly,
as the amount of Gaussian noises $d$ increases keeping other factors fixed,
the mean exit times decreases. As shown in the figure,
the mean exit time is continuous at the boundary
even in the absence of Gaussian noise ($d=0$).
From the numerical simulations (not shown here), we also find that
the mean exit times are continuous at the boundary for $1\leq \a \leq 2$.

\begin{figure}
\begin{center}
\vspace*{-1.0in}
\includegraphics*[width=14.0cm]{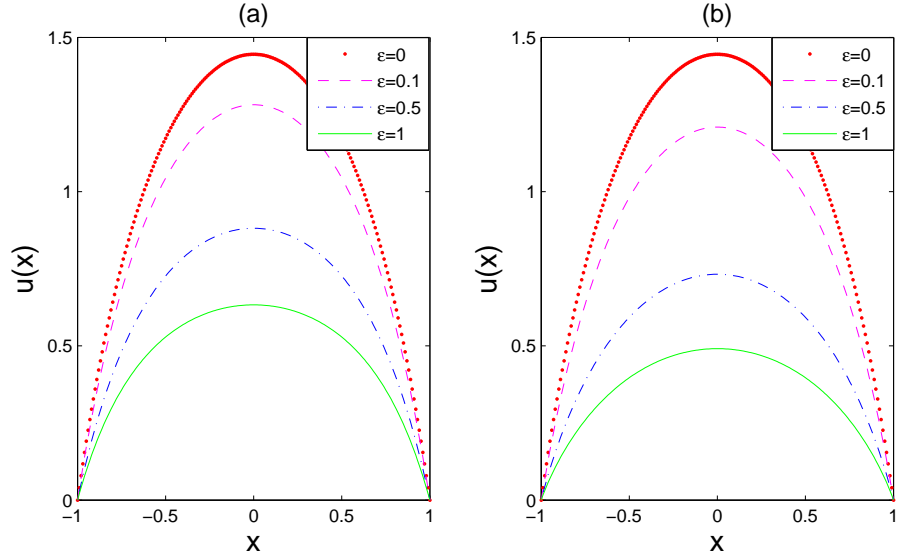}
\end{center}
\caption{Effect of non-Gaussian noise on the mean exit times
of the L\'evy motion defined by the generator in \eqref{AA}.
(a) The mean exit times
for $f(x)=-x, d=1, \a=0.5, D=(-1,1)$ and $\eps=0, 0.1, 0.5, 1$.
(b) Same as (a) except $\a=1.5$.}
\label{engn}
\end{figure}
Next, we look at the effect of non-Gaussian noises by changing the parameter
$\eps$ while keeping $f(x)=-x$, $d=1$ and $D=(-1,1)$. It is obvious
from the numerical results shown in Fig.~\ref{engn} that the mean exit times
decreases when the amount of non-Gaussian noises $\eps$ increases
for all values of $\a$ in $(0,2]$. Due to the presence of significant
Gaussian noises ($d=1$), the shapes of the mean exit times $u(x)$
are parabola-like for all parameter values shown in the figure.
Keeping other parameters fixed, the effect of non-Gaussian noises on mean exit
time is stronger when $\a$ is larger. It is consistent with the previous
result shown in Fig.~\ref{ddsoud1}(a) that,
for small domains, the mean exit times decreases as $\a$
increases.

\subsection{Escape probability}

In this section, we simulate the escape probability described by \eqref{eq.ep}.
In particular, for the special case of $D=(a,b)$
and $E=[b,\infty)$, Eq.~\eqref{eq.ep} becomes
\begin{eqnarray}
  \frac{d}{2} P_E''(x) + f(x) P_E'(x)
  - \frac{\eps C_\a}{\a} \left[\frac{1}{(x-a)^\a}+\frac{1}{(b-x)^\a}\right] P_E(x) &
\nonumber \\
+ \eps C_\a \int_{a-x}^{b-x} \frac{P_E(x+y) - P_E(x) -
   I_{\{|y|<\delta\}} y P_E'(x)}{|y|^{1+\a}}\; {\rm d}y
= -\frac{\eps C_\a}{\a} \frac{1}{(b-x)^\a}, &
\label{eq.ep1}
\end{eqnarray}
for $x \in (a,b)$. The conditions for the escape probability outside
the domain are $P_E(x)=0$ for $x \in (-\infty,a]$
and $P_E(x)=1$ for $x \in [b,\infty)$.

\begin{figure}
\begin{center}
\vspace*{-1.0in}
\includegraphics*[width=14.0cm]{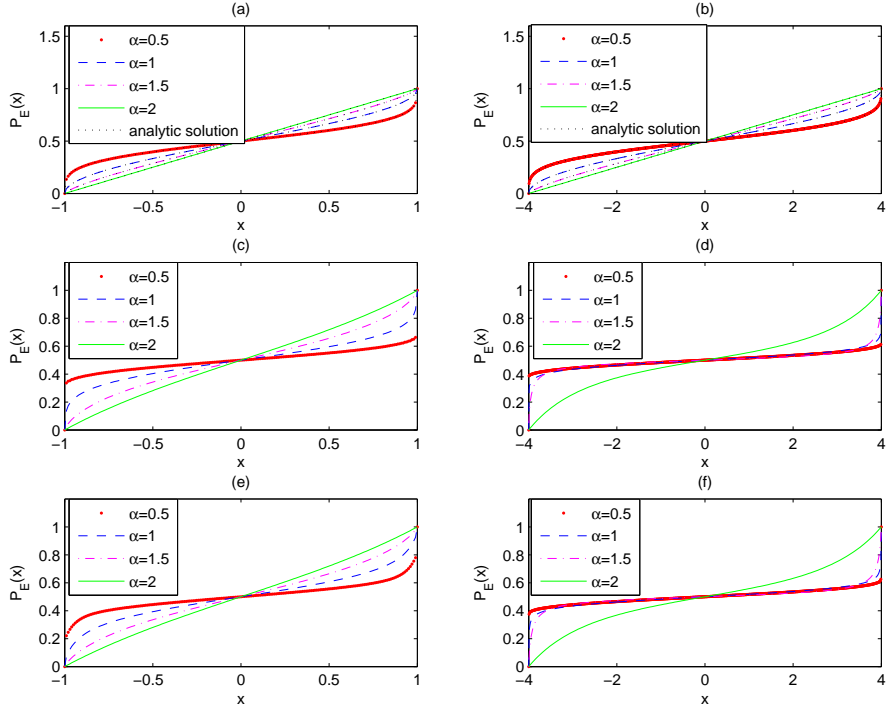}
\end{center}
\caption{The escape probability of the L\'evy motion, defined by the
generator in \eqref{AA}, out of the domain $D$ landing in $E$.
(a) The escape probabilities $P_E(x)$ for the symmetric $\a$-stable
processes, i.e., $d=0$, $f\equiv 0$ and $\eps=1$, with $\a=0.5, 1,
1.5, 2$, $D=(-1,1)$ and $E=[1,\infty)$, and the analytical
results are shown by fine-dotted lines for all $\a$ values;
(b) Same as (a) except
$D=(-4,4)$ and $E=[4,\infty)$; (c) Same as (a) except for the
L\'evy motion with O-U potential $f(x) = -x$, i.e.,
$d=0$, $f(x)=-x$ and $\eps=1$; (d) Same as (c) except $D=(-4,4)$ and
$E=[4,\infty)$; (e) Same as (c) except with added Gaussian
noise, i.e., $d=0.1$, $f(x)=-x$ and $\eps=1$; (f) Same as (e) except
$D=(-4,4)$ and $E=[4,\infty)$. }
\label{ep2}
\end{figure}


First, we verify our numerical schemes by comparing with the analytical
result of the escape probability
for the symmetric $\a$-stable case ($f\equiv 0$, $d=0$, $\eps=1$)
with $D=(-1,1)$ and $E=[1,\infty)$ \cite{Getoor61b}
\begin{equation}
P_E(x) = \frac{(2b)^{1-\a}\Gamma(\alpha)}{[\Gamma(\alpha/2)]^2}
\int_{-b}^x (b^2-y^2)^{\frac{\a}{2}-1} \, {\rm d}y, \quad x\in
(-b,b). \label{eq.asep2}
\end{equation}

As shown in Fig.~\ref{ep2}(a) and (b), the numerical results match
with the analytical results given in Eq.~\eqref{eq.asep2}.
Due to the symmetry of the process and
the domains, the escape probability takes the value of one-half when
the starting point is the position of symmetry $x=0$. The escape
probability is symmetric with respect to the point $(0,1/2)$, i.e.,
\[ P_E(x) + P_E(-x) = 1.\]
These remain true even we add Brownian
noise and the O-U potential to the process.

Due to the symmetry, in the following discussion we focus on positive starting
points in the domain, i.e., $x>0$.
Figure~\ref{ep2} shows that
the probability for the process to escape to the right of the domain
is smaller when the value of $\a$ decreases. For a fixed positive $x$,
the escape probability $P_E(x)$ is the largest in
for the case of Gaussian noise only ($\a=2$), keeping other factors the same.
This property is independent of the domain size or whether there
exists a deterministic driving mechanism $f$.
For larger domain sizes, as shown in Fig.~\ref{ep2}(b,d,f), the escape
probability tends to the value of equal chance $1/2$ especially for small
values of $\a$. By comparing Fig.~\ref{ep2}(c) with Fig.\ref{ep2}(a)
or comparing Fig.~\ref{ep2}(d) with Fig.~\ref{ep2}(b), we find that
the effect of O-U potential
is reducing the escape probability for the same starting point $x$.
The escape probability $P_E(x)$ for the Brownian noise ($\a=2$)
is no longer a straight line in the presence of O-U potential.
Again, for $0<\a<1$, we find that the escape probability
is discontinuous at the boundary of the domain
when the SDE is driven by the O-U potential and "pure" $\a$-stable
symmetric process, as demonstrated by the graphs of $\a=0.5$
in Fig.~\ref{ep2}(c) and (d). Adding Gaussian noises ($d=0.1$) to the processes
increases the chances of escape to the right, as shown in Fig.~\ref{ep2}(e) and (f).

\begin{figure}
\begin{center}
\vspace*{-1.0in}
\includegraphics*[width=14.0cm]{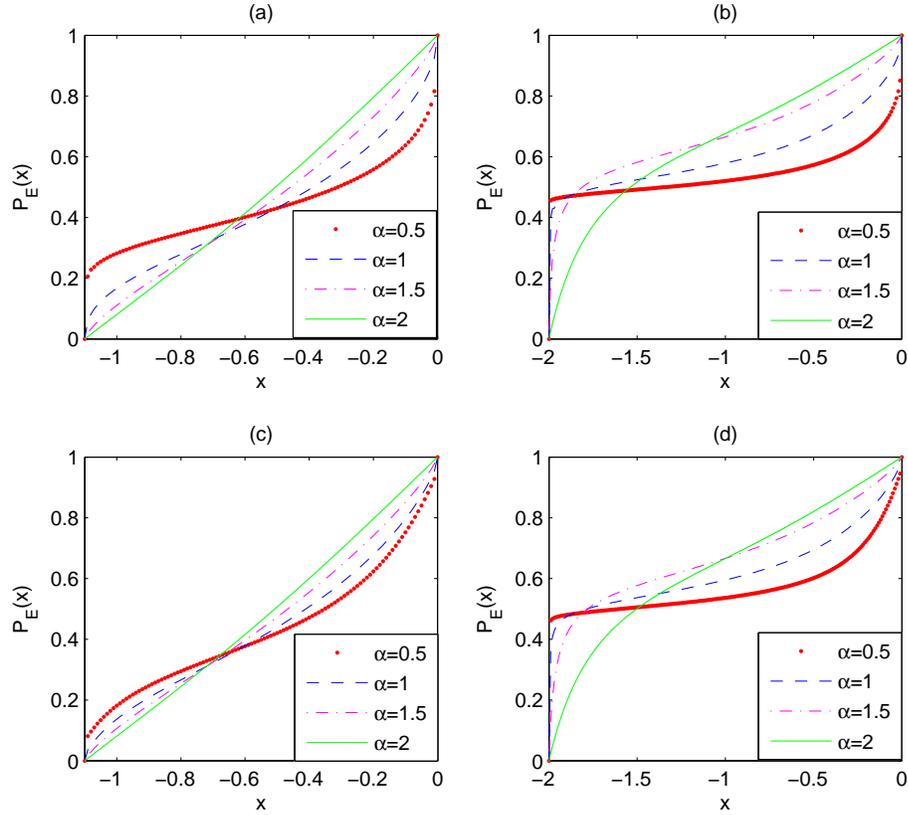}
\end{center}
\caption{The escape probability of the L\'evy motion, defined by the
generator in \eqref{AA}, out of a bounded domain $D$ landing in $E$,
driven by the double-well potential $f(x)=x-x^3$ and noises.
(a) The escape probability $P_E(x)$ for
$d=0, \eps=1$ $D=(-1.1,0)$, $E=[0,\infty)$, $\a=0.5, 1, 1.5, 2$; (b)
Same as in (a) except $D=(-2,0)$; (c) same as (a) except $d=0.1$;
(d) same as (b) except $d=0.1$. }
\label{ep3}
\end{figure}
Next, we consider SDE \eqref{sde} driven by
the double-well potential $f(x)=x-x^3$. The corresponding deterministic
dynamical system ${\rm d}X_t= (X_t-X_t^3)\, {\rm d}t$ has two stable states
located at $x=\pm 1$ while $0$ is an unstable steady state. The double-well
potential is well-known and widely used in phase transition studies.

We investigate the likelihood of the stochastic process
that starts within a bounded domain $D \subset (-\infty, 0)$ and escapes
and lands in the right-half line $E=[0,\infty)$
compared with that of escaping to the left of the bounded domain.
We consider the bounded domains $D$
that includes the left stable point $x=-1$ and the escape-target domain $E$
containing the other stable point $x=1$. When a process lands in $E$,
in absence of the noises, it will be driven to and stay at
the stable point $x=1$. In other words, we try to examine the effect
of noises on the likelihood of the transition from one stable
state to the other.
Figure~\ref{ep3}(a) shows the escape probability
for $D=(-1.1,0)$, $E=[0,\infty)$ and $\a=0.5,1,1.5,2$, when the stochastic
effects are given by $\a$-stable symmetric processes only ($\eps=1,$ and $d=0$).
The escape probability $P_E(x)$ deviates more from a straight line as $\a$
decreases and it is smaller for smaller $\a$ when starting from $x>0.5$.
On the contrast, the probability is larger for smaller $\a$ when the starting
point is close to the left boundary of the bounded domain.
For the bigger domain $D=(-2,0)$ shown in Fig.\ref{ep3}(b),
the likelihood of escape to the right
is more than a half for most of the starting points $x>-1.5$
and $1\leq \a \leq 2$; the probability stays close to a half for most of
the starting points when $\a=0.5$. Note that, in the case of $\a<1$,
the probability is discontinuous at the left boundary of the domain.
As shown in Fig.\ref{ep3}(c) and (d), the differences in escape probabilities
among different values of $\a$ become smaller when an amount of Gaussian noises is added ($d=0.1$), but otherwise probabilities have the similar
values and properties compared with those of $d=0$.

\section{Conclusion}

In summary, we have developed an accurate numerical scheme for solving
the mean first exit time and escape probability for SDEs with non-Gaussian L\'evy motions. We have analyzed
the numerical error due to the singular nature of the L\'evy measure
corresponding to jumps and accordingly, added a correction term to the numerical scheme.
We have validated the numerical method by comparing with analytical
and asymptotic solutions. For arbitrary deterministic driving force,
we have also given asymptotic solutions of the mean exit time
when the pure jump measure in the L\'evy motion is small.

Using both analytical and numerical results, we find that the mean exit
time depends strongly on the domain size and the value of $\a$
in the $\alpha-$stable L\'evy jump measure.
For example, for $\a$-stable L\'evy motion, the mean exit time
can help us decide which of the two competing
factors in $\a$-stable L\'evy motion, the jump frequency or the jump size,
is dominant. The small jumps with high frequency, corresponding
to L\'evy motion with $\a$ closer to $2$, make it easier
to exit the small domains.
On the other hand, it is easier to exit large domains for
$\a$-stable L\'evy motion with $\a$ closer to $0$,
which has the characteristics of the large jumps with low frequency.
Another observation from the results is that, for smaller values
of $\a$ ($0<\a<1$), the mean exit time profiles are flat
away from the boundary of the domain.
This implies that the jump sizes of the processes
are usually larger than the domain sizes, thus the mean exit times
have small variations for different starting positions.

The probability for the process to escape to the right of the domain
is smaller when the value of $\a$ decreases. For a fixed positive $x$,
the escape probability $P_E(x)$ is the largest in
for the case of Gaussian noise only ($\a=2$), keeping other factors the same.
This property is independent of the domain size or whether there
exists a deterministic driving mechanism $f$.
The escape probability is shown to vary significantly
with the underlying vector field.

The mean exit time and escape probability could become discontinuous
at the boundary of the domain, when the process is subject to certain
deterministic potential and the value of $\a$ is in $(0,1)$.



\end{document}